\documentclass[12pt]{article}
\input epsf      
\usepackage{graphicx}
\usepackage{amssymb}
\usepackage{amsfonts}
\usepackage{amsthm}
\usepackage{amsfonts}

\newcommand \half {{\mathbb H}}

\newcommand{\compose}{{\circ}}

\newcommand \ra {\rightarrow}

\newcommand{\ba}[1]{\begin{array}{#1}}
\newcommand{\ea}{\end{array}}

\newcommand{\be}{\begin{equation}}
\newcommand{\ee}{\end{equation}}
\newcommand{\bea}{\begin{eqnarray}}
\newcommand{\eea}{\end{eqnarray}}
\newcommand{\beann}{\begin{eqnarray*}}
\newcommand{\eeann}{\end{eqnarray*}}

\def\reff#1{(\ref{#1})}
\newtheorem{theorem}{Theorem}

\begin{document}

\bibstyle{ams}

\title{
Conformal Invariance and Stochastic Loewner Evolution Predictions 
for the 2D Self-Avoiding Walk - Monte Carlo Tests
}

\author{Tom Kennedy
\\Department of Mathematics
\\University of Arizona
\\Tucson, AZ 85721
\\ email: tgk@math.arizona.edu
}
\maketitle

\begin{abstract}

Simulations of the two-dimensional self-avoiding walk (SAW) are performed in 
a half-plane and a cut-plane (the complex plane with 
the positive real axis removed) using the pivot algorithm. 
We test the conjecture of Lawler, Schramm and Werner
that the scaling limit of the two-dimensional SAW is given by Schramm's 
stochastic Loewner evolution (SLE). The agreement is found to 
be excellent.  
The simulations also test the conformal invariance of the SAW 
since conformal invariance implies that if we map infinite length walks in 
the cut-plane into the half plane 
using the conformal map $z \rightarrow \sqrt{z}$, then the resulting
walks will have the same distribution as the SAW in the half plane. 
The simulations show excellent agreement between the distributions. 

\end{abstract}

\bigskip


\newpage

\section{Introduction}

Lawler, Schramm and Werner  \cite{lsw_saw} have conjectured that 
the scaling limit of the two-dimensional self-avoiding walk (SAW) 
is given by Schramm's \cite{schramm} stochastic Loewner evolution (SLE).
SLE is a two dimensional conformally invariant 
random process which depends on a parameter $\kappa$, and so is 
denoted SLE$_{\kappa}$. 
Chordal SLE refers to the version of SLE in which the random curve or set 
connects two points on the boundary of a simply connected domain. 
It is usually defined first for the case where the domain is the 
half-plane and the two boundary points are $0$ and $\infty$. Its 
definition is them extended to other simply connected domains $D$ and 
boundary points using a conformal map from the half-plane to $D$ 
which maps the two boundary points appropriately. 
If $\kappa < 4$, chordal SLE gives a probability measure on simple curves, 
i.e., curves that do not intersect themselves \cite{rs}.
The conjecture of Lawler, Schramm and Werner is that for any simply 
connected domain $D$ and points $z$ and $w$ on its boundary, SLE$_{8/3}$ is 
the scaling limit of SAW's that go from $z$ to $w$ and stay inside $D$.

For $\kappa=8/3$, Lawler, Schramm and Werner \cite{lsw_cr}
have a theorem that
makes it possible to explicitly compute the distributions of many random
variables associated with the SLE random curve.
For the scaling limit of the SAW, these random variables  can be studied 
by simulation. Thus one can 
numerically test their conjecture that the scaling limit of the 
SAW is SLE$_{8/3}$ by comparing the distributions from simulations
of the SAW with the exact distributions for SLE$_{8/3}$. 
This test was carried out for two such random variables for the SAW in the 
upper half-plane in \cite{tk_saw_sle},
and excellent agreement was found. 
In this paper we consider more random variables for which the 
exact distribution can be computed for SLE$_{8/3}$. 
We compare their exact distributions with the numerical distributions
of the same random variables for the SAW in the half-plane.  
We also simulate the SAW in the cut-plane consisting of the 
complex plane minus the non-negative real axis. 
The map $z \rightarrow \sqrt{z}$ takes the cut-plane onto the half-plane,
and by composing the random variables for the half-plane with 
this map we obtain corresponding random variables for the cut-plane. 
We compare their distributions for the SAW from our simulations for the 
cut-plane with the exact distributions for SLE$_{8/3}$.
We also consider the probability that the walk passes to the right 
of a given point in the half-plane (or the cut-plane) and compare 
this probability for the SAW simulations with an exact formula of Schramm 
\cite{schramm_perc} for SLE.

Note that for both of the domains we consider, the terminal point of
the walk is at infinity. This case is particularly well suited to 
simulations, since it is expected that we can construct the scaling limit
by considering all SAW walks
with a fixed length $N$ which start at the origin, taking the 
limit $N \rightarrow \infty$ and then taking the limit that the 
lattice spacing goes to zero.
For a domain $D$ and two finite points $z$ and $w$ on its boundary,
the scaling limit should be constructed as follows. We introduce 
a lattice and consider all self-avoiding walks 
which start at $z$ and end at $w$. The probability of such a walk
is taken to be proportional to $\beta^{-N}$ where $N$ is the number of steps 
in the walk, and $\beta$ is the constant such that the number of 
SAW's in the plane starting at the origin grows with the number of 
steps, $N$, as $\beta^N$. The measure is normalized so that it is a 
probability measure. We then take the limit of this measure as the lattice
spacing goes to zero. The construction of the scaling limit in the case 
of SAW's with infinite terminal point is rather different from the case
of a finite terminal point, so it would be interesting to test the 
conjecture that the scaling limit is given by SLE$_{8/3}$ in the 
case of a finite terminal point.

In addition to describing the scaling limit of the SAW, 
SLE is conjectured to describe the scaling limit of a large number 
of other two dimensional models. 
Many of these conjectures have been proved recently. 
Schramm showed that if the loop-erased random walk has a conformally invariant
scaling limit, then that limit must be SLE$_2$ \cite{schramm}. 
He also conjectured that the scaling limit of percolation 
should be related to SLE$_6$, and the scaling limit of uniform spanning
trees (UST) is described by SLE$_2$ and SLE$_8$.
The conjectures for the loop-erased random walk and 
the UST have been proved by Lawler, Schramm and Werner \cite{lsw_lerw}.
Smirnov has proved the conformal invariance conjecture
for critical percolation on the triangular lattice and that
SLE$_6$ describes the limit \cite{smirnov}.
Lawler, Schramm and Werner used SLE$_6$ to rigorously determine the 
intersection exponents for Brownian motion 
and proved a conjecture of Mandelbrot that the outer boundary of a Brownian
path has Hausdorff dimension 4/3 \cite{lswa, lswb, lswc}.
The random cluster representation of the Potts model for $0<q<4$ 
was conjectured by Rohde and Schramm to be 
related to the SLE process as well
\cite{rs}. 

\section{SLE predictions}

The random variables we consider are defined for curves in the upper 
half-plane as follows. Note that these random variables are defined 
both for the SAW and for SLE. We use $\gamma$ to denote the 
random curve in both cases.
Consider a horizontal line at a height of $c$ above the horizontal axis. 
The curve $\gamma$ will intersect it, possibly more than once,
and we look for the left-most intersection. The random variable $X_e$ 
is the $x$-coordinate of this intersection, divided by $c$. So 
\be 
X_e={ 1 \over c} \min \{ x : x + ic \in \gamma \}
\ee
We can also consider the first intersection of the curve with the 
horizontal line. (``First'' means the first intersection as we traverse the 
curve starting at the origin.) 
We let $X_f$ be the $x$-coordinate of this 
intersection, divided by $c$. 
(The subscripts $e$ and $f$ are for ``extreme'' and ``first,'' respectively.) 
The next random variable is defined using a vertical line at a 
distance $c$ to the right of the origin. The curve will intersect it, 
and we look for the lowest intersection. The random variable $Y_e$ 
is the $y$-coordinate of this intersection, divided by $c$. So 
\be 
Y_e={ 1 \over c} \min \{ y : c + iy \in \gamma \}
\ee
The random variable $Y_f$ is the $y$-coordinate of the first intersection, 
divided by $c$. 
For the final random variable, consider a semi-circle of radius $c$ 
centered at the the point $cd$ on the real axis where $|d|<1$. So the origin
where the random curve starts is inside the semicircle.   
The intersections of the random curve with the semicircle are of the form 
$c(d+e^{i \theta})$ and we look for the intersection with the smallest 
$\theta$. The random variable $\Theta_e$ is this smallest angle, 
normalized so that it ranges between $0$ and $1$. So 
\be 
\Theta_e ={ 1 \over \pi} \min \{ \theta: c (d+e^{i \theta}) \in \gamma \}
\ee
The random variable $\Theta_f$ is the angle of the first intersection, 
again normalized so that it ranges from $0$ to $1$.
If the probability measure is invariant under dilations, then the 
distributions of all of these random variables are independent of $c$. 
This is true for SLE and is expected to be true for the scaling limit 
of the SAW.

The distributions of $X_e,Y_e,\Theta_e$ are all easily computed using the 
following theorem of Lawler, Schramm and Werner. 
Let $\half$ be the upper half-plane. Let $A$ be a compact subset of the 
closure of $\half$ such that $\half \setminus A$ is simply connected 
and $0$ is not in $A$. 
Let $\Phi_A$ be the conformal map from $\half \setminus A$ onto $\half$ which 
fixes $0$ and $\infty$ and has $\Phi_A^\prime(\infty)=1$.
We continue to denote the random curve generated by SLE, the SLE ``trace,'' 
by $\gamma$. 

\begin{theorem} (Lawler, Schramm, Werner \cite{lsw_cr}.)
For $\kappa=8/3$, chordal SLE in the half plane has 
\be
P(\gamma \cap A = \emptyset) = \Phi^\prime_A(0)^{5/8} 
\ee
\end{theorem} 

Our next step is to use this theorem to compute the distributions of 
$X_e$, $Y_e$ and $\Theta_e$. 

\subsection{Hitting the horizontal line}

It is convenient to take $c=\pi$ to compute the distribution of $X_e$. 
Let $L_t$ be the horizontal ray which starts at $t + \pi i$ and goes to the
left. Let $\Phi_{L_t}(z)$ be the conformal map which maps $\half \setminus L_t$
onto $\half$ and satisfies the conditions in the theorem. 
Note that $X_e \le t / \pi $ if and only if $\gamma$ hits $L_t$. 
So by the theorem 
\be
P(X_e \le t / \pi) = 1 - \Phi_{L_t}^\prime(0)
\label{eqa}
\ee
The map $w(z) = z + \ln(z) +1+t$ maps $\half$ onto $\half \setminus L_t$.
We need the inverse of this map but it cannot 
be explicitly found. The inverse should be normalized so that 
it fixes $0$ and $\infty$ and has derivative $1$ at $\infty$. 
The above map does not fix $0$, but meets the other two conditions.
Fixing $0$ is not necessary since we can achieve this condition
by just adding a constant to the inverse map, and this which will not 
change its derivative. So we have 
\be
\Phi_{L_t}^\prime(0)={dz \over dw}(0) = 
\left({dw \over dz} (z_0) \right)^{-1} 
\ee
where $z_0$ is the image of $0$ under the inverse map, i.e., 
$0= z_0 + \ln(z_0) + 1 + t$.

\begin{figure}
\begin{center}
\includegraphics[width=11cm]{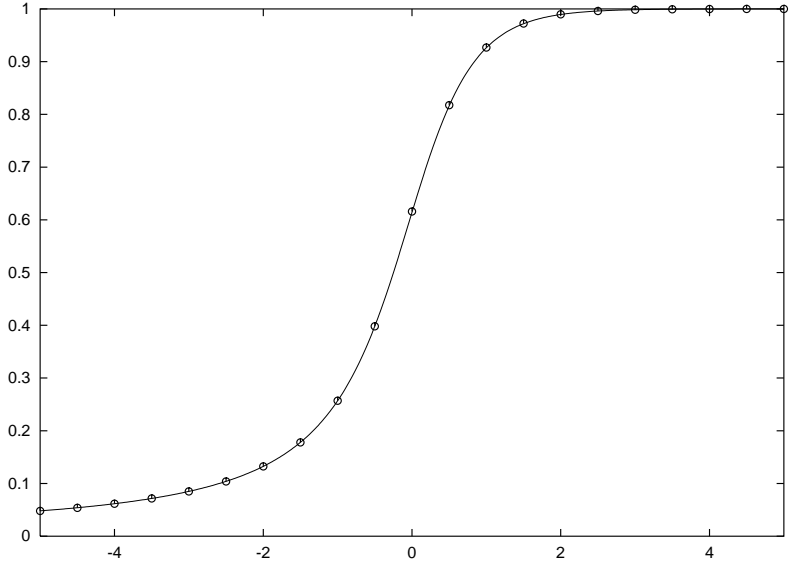}
\caption{ The distribution, $P(X_e \le t)$, of $X_e$ for the half-plane. 
The solid line is the distribution for SLE$_{8/3}$, and the open 
circles are the results of the simulation of the SAW.}
\label{sec_half_x_e}
\end{center}
\end{figure}

Define 
\be
g(x)= x + \ln(x)
\ee
This is an increasing function which maps $(0,\infty)$ onto the real line,
so it has an inverse that maps the real line to $(0,\infty)$.
Note that $z_0=g^{-1}(-t-1)$. We have 
\be
{dw \over dz}(z_0) = 1+ { 1 \over g^{-1}(-t-1)}
\ee
So \reff{eqa} and a trivial change of variables gives 
\be
 P(X_e \le t)
= 1 - \left({g^{-1}(- \pi t - 1) \over g^{-1}(- \pi t - 1) +1}\right)^{5/8}
\label{x_distrib}
\ee
Although $g^{-1}$ cannot be explicitly computed, it can be trivially 
computed numerically. The graph of the above distribution is the solid
line in figure \ref{sec_half_x_e}. The open circles in the figure 
are the results of the simulation for the SAW. 

We can find the asymptotic behavior of the distribution in \reff{x_distrib}
as $t$ goes to $\pm \infty$. 
For large positive $t$,  $g(t)= t + \ln(t) \approx t$.
So as $t \rightarrow -\infty$, $g^{-1}(- \pi t-1) \approx - \pi t$, and so 
\be P(X_e \le t) \approx 1 - \left({- \pi t \over - \pi t +1}\right)^{5/8}
\approx  - {5 \over 8 \pi t }, 
\quad as \quad t \rightarrow -\infty
\ee
As $t \rightarrow 0$, $g(t) \approx \ln(t)$. So as $t \rightarrow \infty$, 
$g^{-1}(- \pi t -1 ) \approx e^{- \pi t - 1}$. So 
\be P(X_e \le t) \approx 1 - \left({e^{- \pi t -1} 
\over e^{- \pi t -1} +1}\right)^{5/8}
\approx 1 - e^{-5(\pi t+1)/8},
\quad as \quad t \rightarrow \infty
\ee
As $t \rightarrow -\infty$, the probability goes to zero slowly, 
but as  $t \rightarrow \infty$, the probability goes to one exponentially
fast. This is reasonable since when $X_e$ is very negative it only means 
there is at least one intersection with the horizontal line far 
to the left of the origin, but when $X_e$ is very positive it means that 
all intersections with the horizontal line are far to the right of the origin.

\subsection{Hitting the vertical line}

\begin{figure}
\begin{center}
\includegraphics[width=11cm]{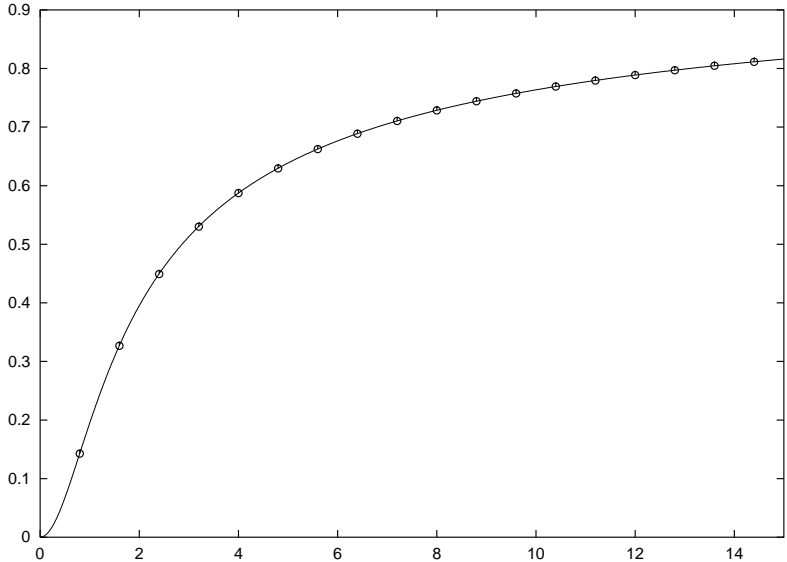}
\caption{The distribution of $Y_e$ for the half-plane. 
The solid line is SLE$_{8/3}$, and the open circles are the SAW.}
\label{sec_half_y_e}
\end{center}
\end{figure}

The distribution of $Y_e$ was studied in \cite{tk_saw_sle}.
We take $A_t$ to be the line
segment from $1$ to $1+it$. The conformal map that maps $\half \setminus A_t$
onto $\half$ with the required normalizations is 
\be 
\Phi_{A_t}(z) = i \sqrt{-(z-1)^2 - t^2}
\ee 
where the square root has its branch cut along the negative real axis. 
Thus the distribution of $Y_e$ is 
\be 
P(Y_e \le t) = P(\gamma[0,\infty) \cap A_t \ne \emptyset) = 
1-\Phi^\prime_{A_t}(0)^{5/8} = 1-(1+t^2)^{-5/16}
\label{y_distrib}
\ee
Figure \ref{sec_half_y_e} shows this distribution and the results of 
the simulation for the SAW. 

\subsection{Hitting the circle}

It is convenient to translate so the semicircle is centered at the origin.
Setting $c=1$, this means the random curves start at $-d$. So
$P(\Theta_e \le t) = 1 - \Phi^\prime(-d)^{5/8}$,
where $\Phi$ is a conformal map which takes the half-plane minus 
the arc $A_\phi=\{e^{i \theta} : 0 \le \theta \le \pi t \}$ onto the 
half-plane with the normalizations that the map fixes $\infty$
and has derivative $1$ at $\infty$. (As in the previous case, we ignore the 
condition that the map fixes the origin since it does not affect the 
derivative.)

The conformal map
\be 
z \ra {z-1 \over z+1}
\ee
sends the upper half-plane (including $\infty$) onto itself,
and it sends the upper half of the unit circle to the upper half of 
the imaginary axis. Let 
\be 
a = {\sin (\pi t) \over 1 + \cos (\pi t)}
\ee
Then the arc $A_t$ is mapped onto the line segment from $0$ to $ia$. 
We can then map $\half$ with this line segment removed onto 
$\half$ as we did in the previous section. Composing these two maps
and multiplying by a factor of $\sqrt{1+a^2}$ for later convenience, 
we define
\be
\psi(z) = i \sqrt{1 + a^2} \left[-{(z-1)^2 \over (z+1)^2}- a^2 \right]^{1/2}
\ee
with the branch cut for the square root being the negative real axis. 
The map $\psi$ sends $\half \setminus A_t$ onto $\half$ .

\begin{figure}
\begin{center}
\includegraphics[width=11cm]{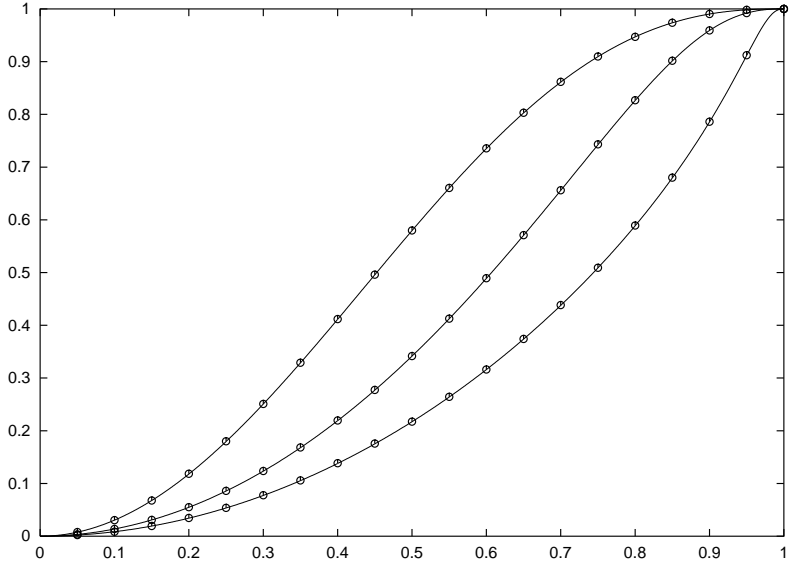}
\caption{The distribution of $\Theta_e$ for the half-plane for $d=0, 0.5, 
0.9$. ($d$ increases from left to right.)  
The solid lines are SLE$_{8/3}$, and the open circles are the SAW.}
\label{sec_half_zzz_e}
\end{center}
\end{figure}

The map $\psi$ does not send $\infty$ to itself.
For $z$ near $\infty$,  
\be
\psi(z) = (1+a^2) \left(1  - {2 \over (1 + a^2) z} + \cdots \right)
\ee
In particular, $\psi(\infty)= (1+a^2)$.
Now let 
\be 
\Phi_{A_t}(z) = {2 \over 1+a^2 - \psi(z)} 
\ee
For large $z$, $\Phi_{A_t}(z) \approx z $, so the derivative at $\infty$ is 
$1$ as required. 

For real $x$ with $-1 < x < 1$, the choice of branch cut leads to 
\be
\psi(x)= - \sqrt{1+a^2} \left[ {(x-1)^2 \over (x+1)^2} + a^2 \right]
\ee
Define
\be
s=(1+a^2)^{-1} = {1 + \cos (\pi t) \over 2}
\ee
This will prove to be a natural variable to use.
We have 
\be
\psi(x)= - {1 \over \sqrt{s}} 
\left[ {(x-1)^2 \over (x+1)^2} + {1 \over s} - 1 \right]^{1/2}
= - {1 \over s} \left[ 1  - {4xs \over (x+1)^2} \right]^{1/2}
\ee
so
\be
\Phi(x)= { 2 \over {1 \over s} - \psi(x) } 
= { 2s \over 1 + \left[ 1  - {4xs \over (x+1)^2} \right]^{1/2}}
\ee
Computing the derivative $\Phi^\prime(-d)$ then yields  
\be 
P(\Theta_e \le t) = 1- 
\left( { 4 s^2 (1+d) \over \left(1-d+[(1-d)^2+4ds]^{1/2}\right)^2
\left((1-d)^2+4ds \right)^{1/2}}\right)^a
\label{theta_d_distrib}
\ee
For $d=0, 0.5$ and $0.9$, this distribution and the results 
of the simulation for the SAW are shown in figure \ref{sec_half_zzz_e}.

\subsection{Passing right}

\begin{figure}
\begin{center}
\includegraphics[width=11cm]{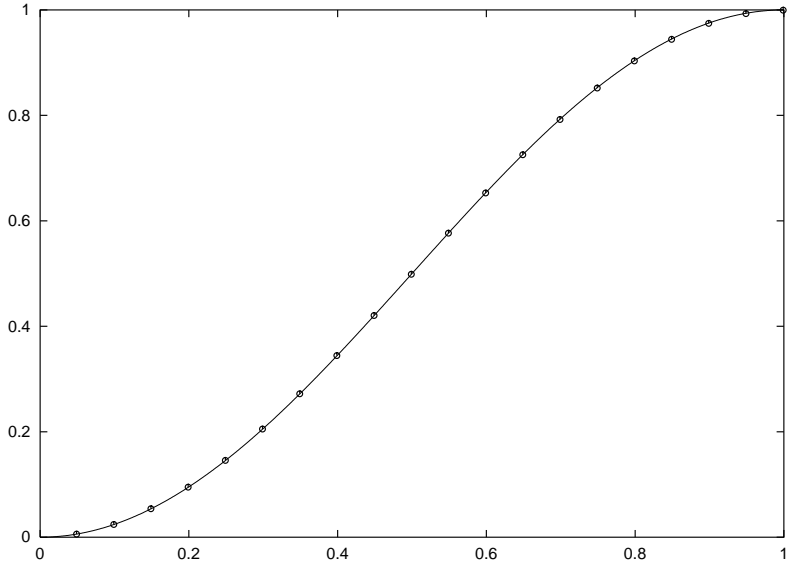}
\caption{The probability that the walk passes to the right of a point 
as function of its polar angle for walks in the half-plane.
The horizontal axis is the angle divided by $\pi$, 
so that it ranges from $0$ to $1$. 
The solid line is the exact result for SLE$_{8/3}$, and the open 
circles are the results of the simulation of the SAW.}
\label{sec_half_pass}
\end{center}
\end{figure}

In addition to the distributions of the 
random variables $X_e$, $Y_e$ and $\Theta_e$, we also
consider the following probability. Fix a point in the 
upper half-plane. One can then ask if the random curve passes to the 
right or left of this point. 
For SLE this probability only depends on the polar
angle of the point since SLE is invariant under dilations. 
This should also be true for the scaling limit of the SAW,
since it is expected to be invariant under dilations.
Schramm \cite{schramm_perc} rigorously derived an explicit formula
for this probability for $\kappa < 8$. 
For general $\kappa$ it is given by a hypergeometric function, but for 
$\kappa=8/3$, his formula is quite simple. Denoting the probability that the 
curve passes to the right of a point with polar angle $\theta$ by $p(\theta)$, 
he showed that for $\kappa=8/3$
\be 
p(\theta) = {1 \over 2} (1-\cos(\theta))
\ee
In our simulations we study this probability by fixing a radius $c$
and computing the probability the path passes to the right of $c e^{i \theta}$
for a large set of values of $\theta$. 
The above function and the results of the SAW simulation are shown
in figure \ref{sec_half_pass}. (Note that the horizontal axis in the 
figure is $\theta/\pi$.)

\subsection{The cut-plane}

The cut-plane we consider is the plane with the non-negative real axis removed.
Let $f(z)=\sqrt{z}$, with the branch cut along the positive real
axis. Then $f$ maps the cut-plane 
onto the upper half-plane. We will continue to 
denote curves in the upper half-plane by $\gamma$, and use $\hat \gamma$
to denote curves in the cut-plane. 
Given a curve $\hat \gamma$ in the cut-plane, 
$\gamma = f \compose \hat \gamma$ 
is a curve in the upper half-plane. 
So we can define the various random variables for the cut-plane 
by applying their definitions in the half-plane to $f \compose \hat \gamma$. 
We will put a $\hat{}$ on top of random variables defined on curves
in the cut-plane. 
For the simulations it is useful to work out these definitions explicitly 
in terms of the curve $\hat \gamma$ in the cut-plane, rather than map 
each SAW in the cut-plane to the half-plane.  

First consider $\hat \Theta_e$ and $\hat \Theta_f$ for $d=0$. 
The map $f$ simply divides the 
polar angle by $2$, so for curves $\hat \gamma$ in the cut-plane,
\be 
\hat \Theta_e 
  ={ 1 \over 2 \pi} \min \{ \theta: c e^{i \theta} \in \hat \gamma \}
\ee
The random variable $\hat \Theta_f$ is the polar angle of the 
first intersection of $\hat \gamma$ with the circle, divided by $2 \pi$. 
If $d \ne 0$, the image of the semicircle under $z \ra z^2$ is not a 
circle. We have not simulated $\hat \Theta_e$  or $\hat \Theta_f$  
in this case. 

To find the definition of $\hat X_e$, we first take $c=1$.
The image of the horizontal line $\{i+t: -\infty < t < \infty \}$ 
under $z \rightarrow z^2$ is a parabola whose axis is the horizontal axis
and which opens to the right, 
\be
x=t^2-1, \quad  \quad  \quad y=2t
\ee 
In the half-plane, $X_e$ is the smallest $t$ such that $i+t$ is on the curve.
In the cut-plane, we consider all intersections of $\hat \gamma$ with 
the parabola and find the intersection with the smallest $y$-coordinate.
Since $t=y/2$, $\hat X_e$ is one half of the $y$-coordinate of this 
``lowest'' intersection. Equivalently,
\be
\hat X_e = \min \{ t : (t^2-1,2t) \in \hat \gamma \} 
\ee
SLE is invariant under dilations of the cut-plane, and the scaling limit
of the SAW in the cut-plane is expected to have this invariance as well.
So for $c \ne 1$ we can take the parabola to be
\be
x=c(t^2-1), \quad  \quad  \quad y=2ct
\ee 
and let 
\be
\hat X_e = \min \{ t : (c(t^2-1),2ct) \in \hat \gamma \} 
\ee
$\hat X_f$ is the $y$-coordinate of the first intersection
of $\hat \gamma$ with the parabola divided by $2c$.   

To find the definition of $\hat Y_e$, we consider the image of 
$\{1+it: 0 < t < \infty \}$ under $z \ra z^2$. It is the upper half of
a parabola whose axis is the horizontal axis and which opens to the 
left: 
\be
x=1-t^2, \quad  \quad  \quad y=2t, \quad \quad t >0
\ee 
In the half-plane, $Y_e$ is the smallest $t$ such that $1+it \in \gamma$,
so in the cut-plane $\hat Y_e$ is one half of the $y$-coordinate of 
the lowest intersection of $\hat \gamma$ and the half parabola. 
\be
\hat Y_e = \min \{ t : (1-t^2,2t) \in \hat \gamma, t > 0 \} 
\ee
More generally, we can let
\be
\hat Y_e = \min \{ t : (c(1-t^2),2ct) \in \hat \gamma, t > 0 \} 
\ee
$\hat Y_f$ is the $y$-coordinate of the first intersection with the parabola
divided by $2c$.   

We have defined the random variables in the cut-plane so 
that if the probability measure is conformally invariant, 
then they will have the same distribution as their 
counterparts in the half-plane.
Rather than compare the distributions of the 
random variables $X_e, Y_e$ and $\Theta_e$
with those of $\hat X_e, \hat Y_e$ and $\hat \Theta_e$, we will compare 
all these distributions with the SLE$_{8/3}$ predictions, 
eqs. \reff{x_distrib}, \reff{y_distrib} and \reff{theta_d_distrib}.
This tests both the conjecture that the scaling limit of the SAW is
SLE$_{8/3}$ and the conformal invariance of the SAW. 
For the random variables $X_f, Y_f$ and $\Theta_f$, we do not know their
distributions for SLE$_{8/3}$. So we will directly compare the 
distributions of $X_f, Y_f$ and $\Theta_f$ with those of 
$\hat X_f, \hat Y_f$ and $\hat \Theta_f$. This tests the conformal invariance 
of the SAW. 

\section{The simulations}

In all of our simulations the walks had one million steps. 
For the half-plane we ran the pivot algorithm for 10 billion iterations
of the Markov chain. For the cut-plane we ran for 11.4 billion iterations.
The simulation of $\Theta_e$ for $d \ne 0$ in the half-plane was done
separately and consisted of 6.8 billion  iterations. For walks with a million 
steps only about $5\%$ of the proposed pivots are accepted. 
Of course, accepted pivots do not produce independent walks and for 
the random variables considered here most accepted pivots do not even 
change the value of the random variables. So the number of effectively
independent samples is considerably less than the number of accepted pivots. 
Each of the simulations requires about a month on a 1.5 GHz PC.
The exact speed of the simulation depends on the choice of the half-plane vs. 
cut-plane and how many random variables are simulated. 

A walk with $N$ steps is typically of size $N^{3/4}$, so to study the 
various random variables we take $c=l N^{3/4}$, where $l$ is fairly small. 
Note that if we rescaled to make $c$ equal to 1, the lattice spacing
would be $(l N^{3/4})^{-1}$. We will refer to this quantity as the 
``effective lattice spacing.'' 
Note that $l$ is the ratio of the scale used to define the random variable 
to the scale of the walk. 
So we must take $l$ small to make the effect of the finite length of our 
walks negligible. But as $l$ gets smaller, the effective lattice 
spacing gets larger. There is a second effect as $l$ gets smaller. 
For smaller $l$, the fraction of the pivots that change the values of 
the random variables is smaller. So the statistical errors get larger
as $l$ gets smaller.
We do not know a priori what value of $l$ will 
be optimal, so we compute the distributions of each random variable
for four different values of  $l$ in our simulations. 
The particular values of $l$ that we use are determined by some 
experimentation  with much shorter simulation runs. 
We do not use the same four values of $l$ for the different random 
variables.

\begin{figure}[pt]
\begin{center}
\includegraphics[width=11cm]{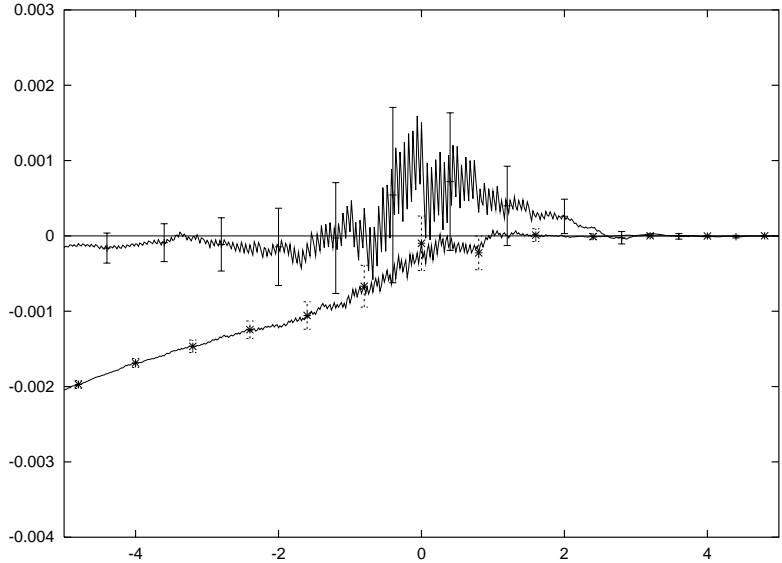}
\caption{Half-plane: The distribution of $X_e$ for the SAW minus the 
distribution of $X_e$ for SLE$_{8/3}$. 
The top curve, with the larger error bars
drawn with solid lines, has $l=0.01$, and the bottom curve has $l=0.05$.
}
\label{sec_half_x_e_dif}
\end{center}
\end{figure}

\begin{figure}[pt]
\begin{center}
\includegraphics[width=11cm]{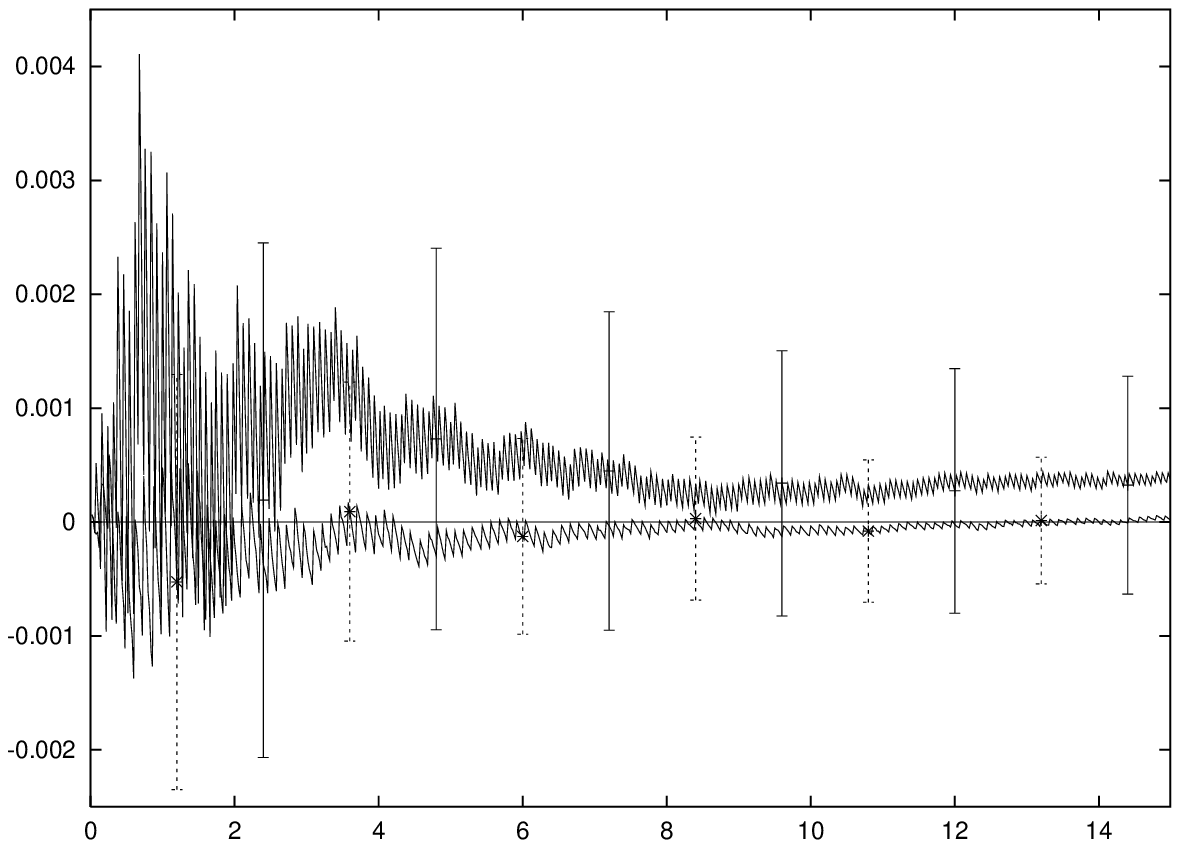}
\caption{Half-plane: The distribution of $Y_e$ for the SAW minus the 
distribution of $Y_e$ for SLE$_{8/3}$. 
The top curve, with the larger error bars
drawn with solid lines, has $l=0.002$, and the bottom curve has $l=0.005$.
}
\label{sec_half_y_e_dif}
\end{center}
\end{figure}

\begin{figure}[pt]
\begin{center}
\includegraphics[width=11cm]{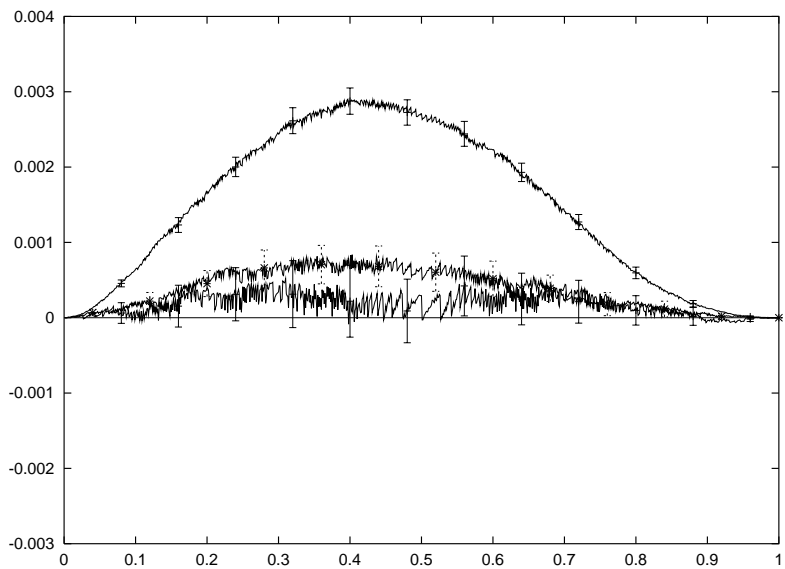}
\caption{Half-plane: For $d=0$, the distribution of $\Theta_e$ 
for the SAW minus the distribution of $\Theta_e$ for SLE$_{8/3}$. 
The three curves shown are for $l=0.2,0.1,0.05$ (from top to bottom).
As $l$ decreases the finite length effects decrease, but the error bars
and lattice effects grow larger.
}
\label{sec_half_z_e_dif}
\end{center}
\end{figure}

\begin{figure}[pt]
\begin{center}
\includegraphics[width=11cm]{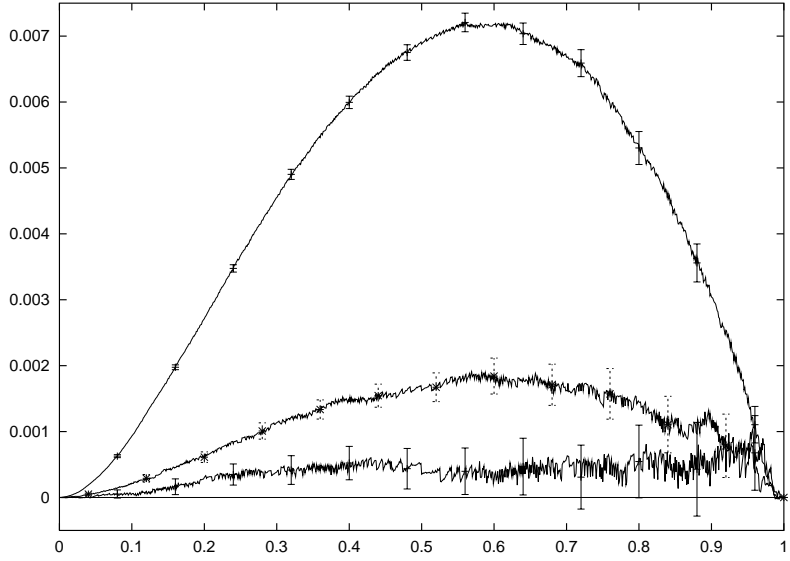} 
\caption{Half-plane: For $d=0.9$, the distribution of $\Theta_e$ 
for the SAW minus the distribution of $\Theta_e$ for SLE$_{8/3}$. 
The three curves shown are for $l=0.2,0.1,0.05$, in order from top to bottom.
}
\label{sec_half_z_e3_dif}
\end{center}
\end{figure}

\begin{figure}[pt]
\begin{center}
\includegraphics[width=11cm]{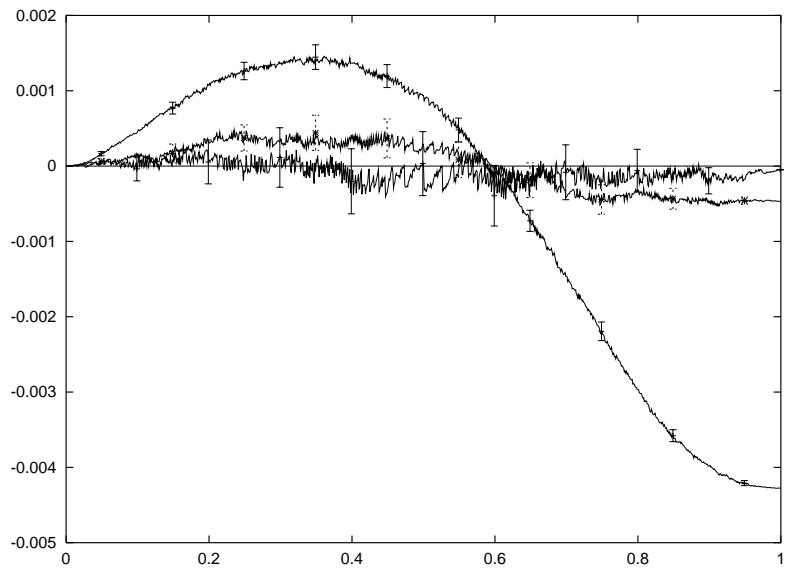}
\caption{Half-plane:The probability that the SAW  passes to
the right of a point as function of the polar angle of the point.
The corresponding function for SLE$_{8/3}$ has been subtracted off.
Going from top to bottom on the left half of the figure, the curves
are $l=0.2,0.1,0.05$. 
}
\label{sec_half_pass_dif}
\end{center}
\end{figure}

In figures \ref{sec_half_x_e} to \ref{sec_half_zzz_e}
we show the distributions of $X_e$, $Y_e$ and $\Theta_e$. 
(Throughout this paper we work with the cumulative distributions of 
our random variables rather than their densities since any simulation 
computes cumulative distributions. Computing 
densities requires taking numerical derivatives of the cumulative 
distributions, and so the densities would have larger statistical errors.)
The solid curves are the exact distributions for SLE$_{8/3}$. 
The circles are the results of the simulation of the SAW.
Figure \ref{sec_half_pass}
studies the probability that the walk passes to the right of 
a point in the upper half-plane as a function of the polar angle
of the point. The solid curve is Schramm's exact result for SLE$_{8/3}$, 
and the circles are the results of the SAW simulation. 
In all of figures  \ref{sec_half_x_e} to \ref{sec_half_pass},
one cannot see any difference between 
the SAW simulations and the exact curves for SLE$_{8/3}$. In figures
\ref{sec_half_x_e_dif} to \ref{sec_half_pass_dif}
we plot the same four quantities, except that now we plot the 
result of the SAW simulation minus the SLE$_{8/3}$ 
functions. The first thing that should be observed in these figures
is the scale of the vertical axis. It is quite small. 
In all but one of these figures the total vertical range shown is 
$0.007$ or $0.7\%$. In figure \ref{sec_half_z_e3_dif} it is $0.008$.

In figures \ref{sec_half_x_e_dif} to \ref{sec_half_pass_dif}
several values of $l$ are shown. 
The nonzero effective lattice spacing means that we are simulating discrete
random variables. So their distributions will be discontinuous. After 
subtracting off the continuous SLE distribution, the jumps will appear in 
the difference as rapid oscillations.
As $l$ increases, the effective lattice spacing decreases, 
and so the oscillations are usually ``faster'' but smaller in amplitude.
Also, for larger $l$ a larger fraction of the pivots change the values of 
the random variables, and so a larger $l$ typically produces smaller 
statistical errors. Both of these effects can be seen in all four of the 
plots.  

\clearpage

\clearpage

\clearpage

As $l$ becomes larger, the effect of the 
finite length of the walk will begin to be seen. 
This is well illustrated by 
figure \ref{sec_half_z_e_dif}, which shows the distribution of $\Theta_e$  
for $d=0$. 
For the largest value of $l$ shown, $l=0.2$, the effect of the finite
length of the walk is clear - the curve differs from zero 
by several times the size of the statistical errors. 
This curve is the smoothest of the three curves and has the 
smallest statistical errors. For $l=0.1$ the finite length effect is
greatly reduced, but is still statistically significant. The $l=0.05$
curve seems to be the best of the values of $l$ that were simulated. 
The maximum difference of the SAW and SLE$_{8/3}$ distributions is only 
about $0.05\%$. Our simulations included a fourth value of $l$ which is 
not shown, $l=0.02$. For this value the larger effective lattice spacing 
and larger statistical errors produce a difference curve that is rougher 
and larger than the $l=0.05$ curve. 
The behavior in figure  \ref{sec_half_z_e3_dif} for the distribution of 
$\Theta_e$ for $d=0.9$ is quite similar to figure \ref{sec_half_z_e_dif},
except that the nonzero $l$ effects appear to  be larger.

In figure \ref{sec_half_x_e_dif} the finite length effect is clearly seen 
in the $l=0.05$ curve; for large negative values of $t$ the deviation of 
this curve from zero is caused by the walk being too short. 
In figure \ref{sec_half_y_e_dif} there are no obvious finite length 
effects; the deviation of the curve from zero
appears to be primarily caused by the nonzero effective lattice spacing.
The deviation is of the same order as the error bars and the oscillations. 
In figure \ref{sec_half_pass_dif} the finite length effects and 
nonzero effective lattice spacing effects are similar to those
seen in figure \ref{sec_half_z_e_dif}. Note that the $l=0.2$ and $l=0.1$
curves are significantly different from zero at the right, corresponding 
to a polar angle of $\pi$. This effect is a result of the nonzero
probability that the walk does not reach the semi-circle or that it 
crosses it, but ends inside the semi-circle. In both of these cases it is 
unclear whether the walk will pass to the right or left of the points on 
the semicircle. The algorithm must make some arbitrary choices in these cases.

\begin{figure}
\begin{center}
\includegraphics[width=11cm]{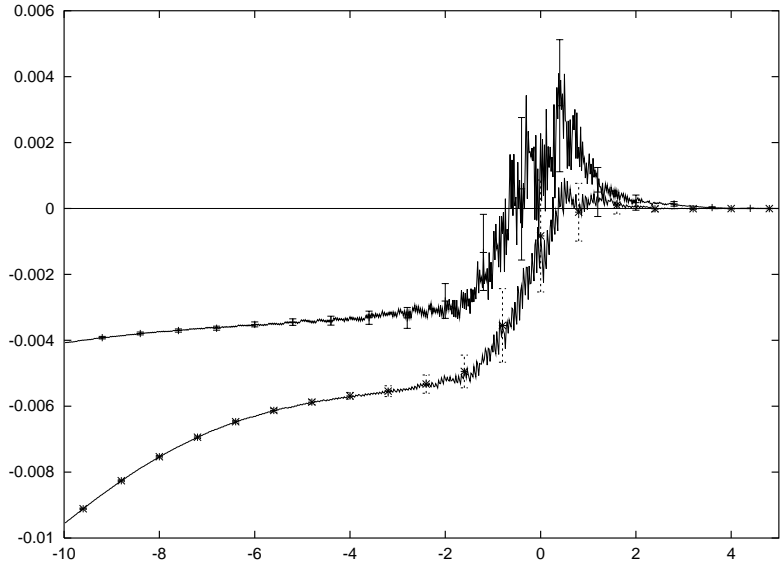}
\caption{Cut-plane: The distribution of $\hat X_e$ for the SAW minus the 
distribution of $\hat X_e$ for SLE$_{8/3}$. 
The top curve, with the larger error bars
drawn with solid lines, has $l=0.002$, and the bottom curve has $l=0.005$.
}
\label{sec_cut_x_e_dif}
\end{center}
\end{figure}

\begin{figure}
\begin{center}
\includegraphics[width=11cm]{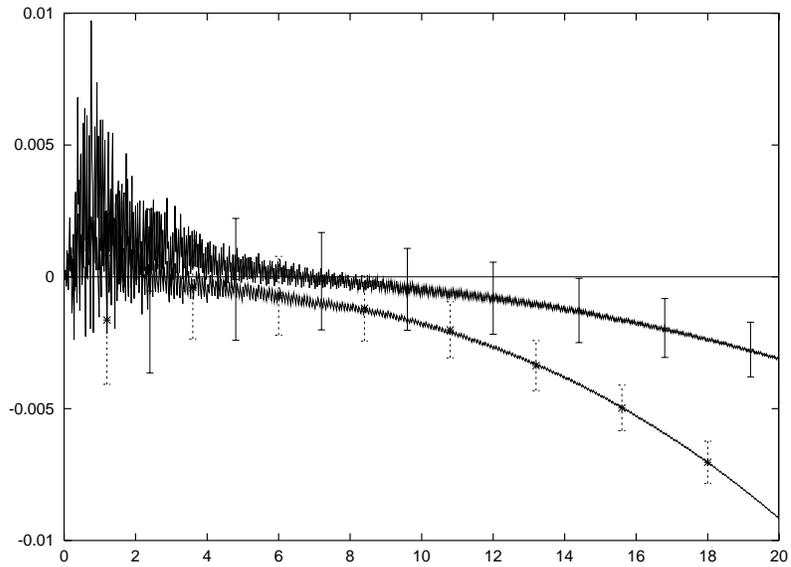}
\caption{Cut-plane: The distribution of $\hat Y_e$ for the SAW minus the 
distribution of $\hat Y_e$ for SLE$_{8/3}$. 
The top curve, with the larger error bars, has $l=0.0005$, 
and the bottom curve has $l=0.001$.
}
\label{sec_cut_y_e_dif}
\end{center}
\end{figure}

\begin{figure}
\begin{center}
\includegraphics[width=11cm]{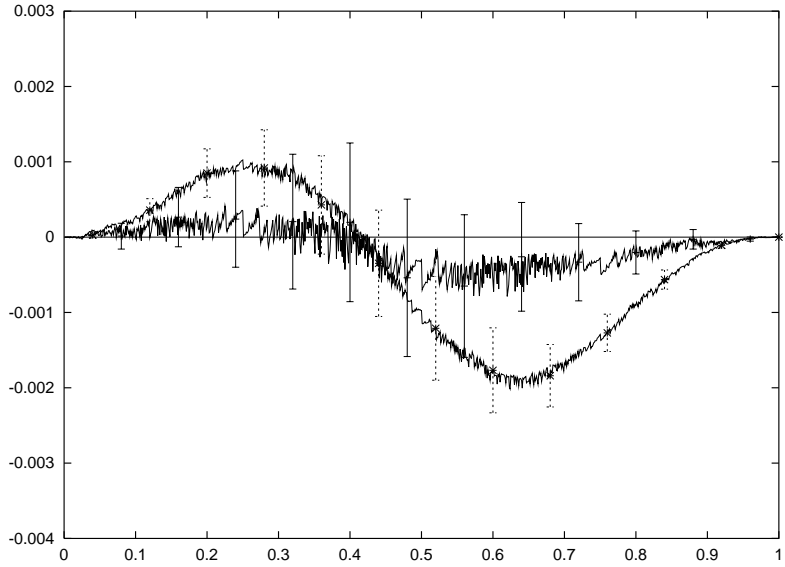}
\caption{Cut-plane: The distribution of $\hat \Theta_e$ for the SAW minus the 
distribution of $\hat \Theta_e$ for SLE$_{8/3}$. 
The curve with the greater deviation from the horizontal axis and 
the error bars drawn with dashed lines has $l=0.05$. The other curve
has $l=0.02$. 
}
\label{sec_cut_z_e_dif}
\end{center}
\end{figure}

\begin{figure}
\begin{center}
\includegraphics[width=11cm]{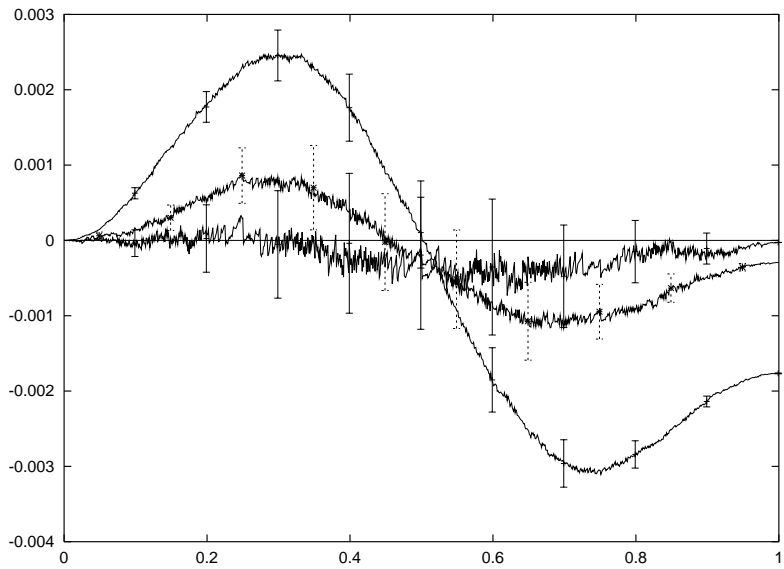}
\caption{Cut-plane: The probability that the walk passes to
the right of a point as function of the polar angle of the point.
Going from top to bottom on the left half of the figure, the curves
have $l=0.1,0.05,0.02$.
}
\label{sec_cut_pass_dif}
\end{center}
\end{figure}

Figures \ref{sec_cut_x_e_dif} through \ref{sec_cut_pass_dif}
show the same quantities as figures
\ref{sec_half_x_e_dif} to \ref{sec_half_z_e_dif}
and \ref{sec_half_pass_dif}, but for
the cut-plane. For the random variable $\hat \Theta_e$ 
(figure \ref{sec_cut_z_e_dif}) and the 
probability of passing right of a point (figure \ref{sec_cut_pass_dif}), 
the agreement is again excellent. 
In both of these figures the vertical scale is $0.007$, the
same as in the corresponding figures for the half-plane. 
For the random variables $\hat X_e$ and 
$\hat Y_e$, figures \ref{sec_cut_x_e_dif} and \ref{sec_cut_y_e_dif},
the agreement is not quite as good, but 
the deviations from the SLE results are still small. 
(In these two figures the vertical scale is two to three times larger
than in the other figures.) 
For these two random
variables it is harder to do accurate simulations for the following 
reason. In the cut-plane, $\hat X_e$ and $\hat Y_e$ depend on the intersections
of the random curve with parabolas. It typically takes a longer length of curve
to attain these intersections than for the lines involved in the 
definition of $X_e$ and $Y_e$ in the half-plane. So in the cut-plane we
must use smaller values of $l$. For $\hat X_e$ in the cut-plane, 
the curves shown use $l=0.002$ and $l=0.005$ as compared to $l=0.01$ and
$l=0.05$ for $X_e$ in the half-plane.
Even with these small values of $l$, the finite length
effects are still quite visible in figure \ref{sec_cut_x_e_dif}.
The deviation of the 
curves from $0$ for the most negative values of $t$ is pronounced.
This is the part of the distribution that is particularly sensitive to 
the need for very long walks to hit the parabola. 
Of course, small values of $l$ mean a large effective lattice spacing 
and large statistical errors. 
For $\hat Y_e$ in the cut-plane, the values of $l$ shown are $0.0005$
and $0.001$, as compared to $0.002$ and $0.005$ for $Y_e$ in the 
half-plane. 
The finite length effects in figure \ref{sec_cut_y_e_dif}
can be seen in the substantial
deviation of the curves from $0$ for large $t$, again a reflection of 
the need for long walks to reach the parabola.  

\clearpage 

\clearpage 

\clearpage 

The scaling limits for the SAW in the half and cut-planes are conjectured 
to be related by the conformal transformation, but there is no reason
that the finite length effects in the two cases should be related. 
Indeed, the simulations show they are quite different. For example, 
compare the curves for the largest values of $l$ in figures 
\ref{sec_half_z_e_dif} and \ref{sec_cut_z_e_dif}.
The curve in figure \ref{sec_half_z_e_dif}  is always positive, 
looking roughly like the first half of a sine wave, 
while the curve in figure \ref{sec_cut_z_e_dif} is both positive and negative. 

Finally, we consider the random variables $X_f$, $Y_f$ and $\Theta_f$
in the half and cut-planes.  
We don't know the exact distributions of these random variable
for SLE$_{8/3}$, but we can still compare the distributions we get from
the simulations of the SAW in the half-plane with the simulations
for the cut-plane. Recall that $\hat X_f$, $\hat Y_f$ and $\hat \Theta_f$ 
(the random variables in the cut-plane) were defined so that they will 
have the same distribution as their counterparts in the half-plane
if the SAW is conformally invariant. 
If we simply plot the distributions themselves,
they agree so well that the difference cannot be seen in the plots.
So instead of plotting the distributions, we plot the distributions
minus various reference functions. These reference functions 
are quite ad hoc. They are chosen to be simple functions that 
are relatively good approximations to the distributions. 
They are defined as follows. 
For $X_f$ and $\hat X_f$ we use the function
\be
F(t)=  {1 \over 2} \left( \tanh(1.16 t) + 1 \right) 
\label{ref_func_x} 
\ee
For $Y_f$ and $\hat Y_f$ we use the distribution of $Y_e$ for 
SLE$_{8/3}$, i.e., 
\be
F(t)= 1-(1+t^2)^{-5/16}
\label{ref_func_y} 
\ee
For $\Theta_f$ and $\hat \Theta_f$ we use 
\be
F(t)=  t -0.12 \sin(2 \pi t) -0.009 \sin(4 \pi t) 
\label{ref_func_theta} 
\ee
We emphasize that these are not meant to be highly accurate 
approximations of the distributions of $X_f$, $Y_f$ and $\Theta_f$. 
One could find better approximations with more complicated functions.
The only purpose of these functions is to provide a convenient 
reference with respect to which we can plot the distributions for 
the half and cut-planes and compare them.

\begin{figure}[tbp]
\begin{center}
\includegraphics[width=11cm]{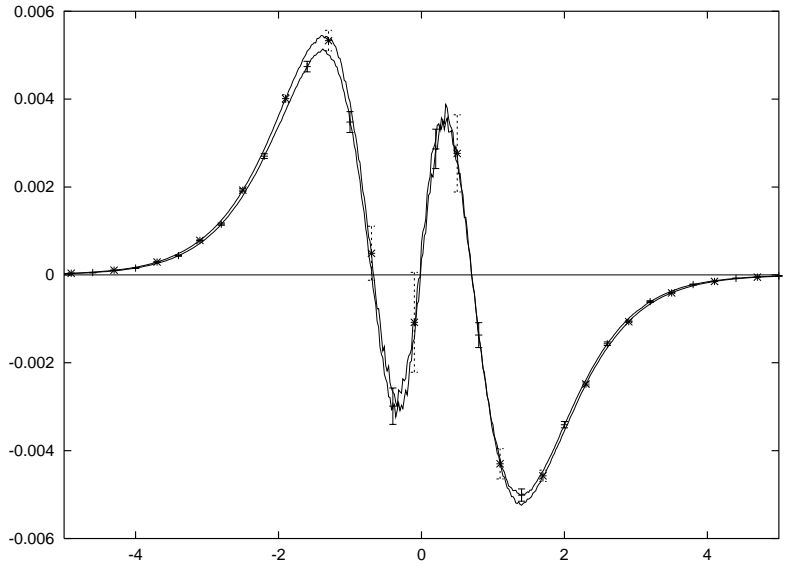}
\caption{For the half and cut planes the distribution of $X_f$ for the 
SAW simulation minus the reference function \reff{ref_func_x} is shown.
In the half-plane $l=0.05$, and in the cut-plane $l=0.02$. 
}
\label{sec_both_x_f_dif}
\end{center}
\end{figure}

\begin{figure}
\begin{center}
\includegraphics[width=11cm]{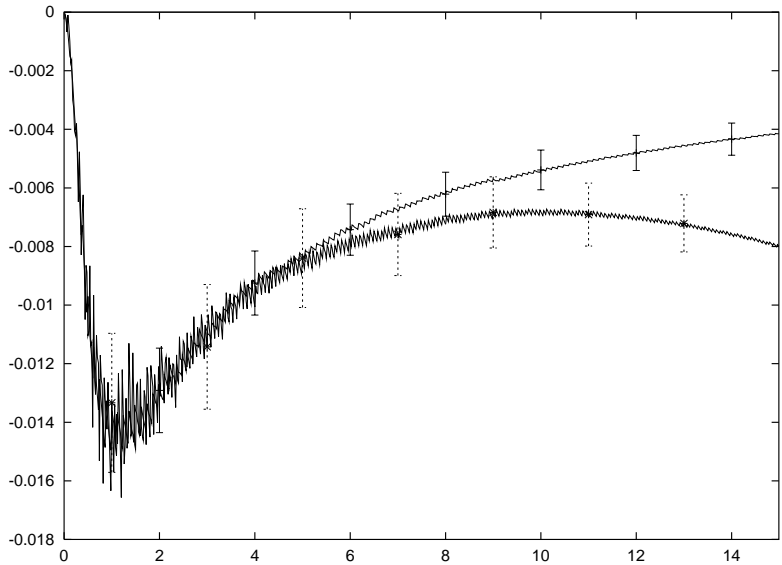}
\caption{For the half and cut planes the distribution of $Y_f$ for the 
SAW simulation minus the reference function \reff{ref_func_y} is shown.
In the half-plane simulation $l=0.005$. In the cut-plane simulation 
$l=0.001$. 
Even with these small values of $l$, the finite length 
effects produce a noticeable difference between the curves for large $t$.
}
\label{sec_both_y_f_dif}
\end{center}
\end{figure}

\begin{figure}
\begin{center}
\includegraphics[width=11cm]{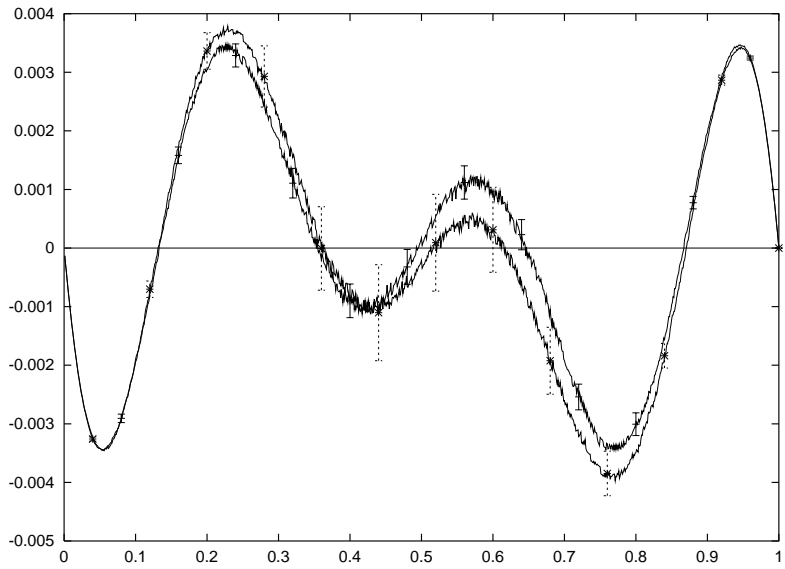}
\caption{For the half and cut planes the distribution of $\Theta_f$ for the 
SAW simulation minus the reference function \reff{ref_func_theta} is shown.
The half-plane simulation used $l=0.1$ and the cut-plane simulation 
used $l=0.05$. The half-plane curve has error bars drawn with solid lines,
while the cut-plane uses dashed error bars.
}
\label{sec_both_z_f_dif}
\end{center}
\end{figure}

Figures \ref{sec_both_x_f_dif} to \ref{sec_both_z_f_dif}
compare the distributions of $X_f$, $Y_f$ and $\Theta_f$ in the 
half-plane with their analogs for the cut-plane. Again, the most
important features of these graphs is the small scale of the vertical
axis. For $X_f$ and $\Theta_f$ the difference between the distributions
in the half and cut-planes is very small. For $Y_f$ the difference is 
somewhat larger for large values of $t$, but still small. We attribute
this greater difference to the larger finite length effects in the 
cut-plane. It can take a walk in the cut-plane a long time to 
reach the parabola involved in the definition of $\hat Y_f$.

\clearpage

\clearpage

\section{Algorithmic considerations}

The pivot algorithm is used for our simulations. (This algorithm is 
discussed in \cite{ms}.) The algorithm picks
a site at random along the walk, called the pivot point,
and picks a random element of the group of symmetries of the lattice 
about the point. This group element is applied to the part of the 
walk after the pivot point. The result is a new nearest neighbor walk,
but it need not be self-avoiding or lie in the upper half-plane. 
The walk is accepted only if both of these conditions are meet. 
Otherwise the proposed walk is rejected and the current walk is counted as
another state in the Markov chain.  The Markov chain trivially 
satisfies detailed balance. In the appendix we show that it is irreducible. 

The speed of the pivot algorithm is typically measured by considering
the average time needed to produce an accepted pivot. 
The algorithm may be implemented \cite{tk_pivot}
so that this time grows with the 
number of steps, $N$, as $O(N^q)$ with $q<1$. 
The exact value of $q$ is not known and probably depends on details of 
the implementation, but simulations indicate the implementation in 
\cite{tk_pivot} has $q<0.57$ in two dimensions.
(This estimate is based on simulations of the walk in the full plane, 
not the half or cut planes.) 

There are two main steps in the pivot algorithm, and both 
would seem to require a time $O(N)$ per accepted pivot. 
The first is the test for self intersections
to see if the new walk should be accepted. The second is actually 
carrying out the pivot. To test for self-intersections 
quickly, we take advantage of
the fact that the walk $\omega$ only takes nearest neighbor steps.
Rather than simply checking if $\omega(i)=\omega(j)$, 
we compute the distance $d=||\omega(i)-\omega(j)||_1$. 
If $d$ is nonzero then we can 
conclude not just that $\omega(i) \ne \omega(j)$, but also that 
\be
\omega(i^\prime) \ne \omega(j^\prime), \quad if \quad |i-i^\prime| 
+ |j-j^\prime| < d
\label{observe}
\ee
Thus we can rule out a large number of potential self intersections 
if $d$ is large.
Since it takes a time $O(N)$ to simply write down a walk with $N$ steps,
the second step of carrying out the pivot would seem to require 
a time that is $O(N)$ per accepted pivot.  
To do better, the key idea is to not 
carry out the pivot each time a pivot is accepted. Instead we
keep track of which pivots have been accepted and only carry them 
out after a certain number have been accepted. 
Details of this implementation of the pivot algorithm may be found
in \cite{tk_pivot}.


In the usual implementation of the pivot algorithm
one chooses the pivot point by giving equal probability to all the 
points on the walk. One can, however, take the probability of picking the 
$i$th site along the walk to be $p(i)$, where $p(i)$ is a function whose 
sum is $1$. The only constraint is that the $p(i)$ must be positive.
If one is interested in the distribution of the end-point of the walk,
then every accepted pivot changes this random variable. 
For this random variable it does not appear that anything could be gained by 
making $p(i)$ non-uniform. However, there is a substantial benefit
to using a non-uniform $p(i)$ for the random variables in this paper. 
All of our random variables typically depend only on a short
segment of the walk near the origin. (The smaller $l$ is, the shorter the 
segment.) So most accepted pivots do not produce any change in the 
random variable. This suggests that it might be worthwhile to choose
pivot locations near the start of the walk more often than pivot locations far
from the start. For the simulations in this paper we define $p(i)$ 
as follows
\be
p(i) = \cases{8c,&if $ 0 \le i < {1 \over 5} N$ \cr
                   4c,&if $ {1 \over 5} N \le i < {2 \over 5} N$ \cr
                   2c,&if $ {2 \over 5} N \le i < {3 \over 5} N$ \cr
                    c,&if $ {3 \over 5} N \le i < N$ \cr}
\ee
where $c={5 \over 16} N^{-1}$ so that the sum of the $p(i)$ is $1$.
This is a rather ad hoc choice, but a crude test indicates that 
for a given number of iterations of the algorithm, it typically reduces the
standard deviation of the random variable by a factor of two. 
A systematic study of the effect of $p(i)$ would be useful.

For each of the six random variables we consider four different values
of $l$. We also consider four values of $l$ for the probability 
of passing to the right of a given point. Thus there are 28 different 
observables to be computed, and some care is necessary to be sure that
the time required for this part of the simulation does not dominate the 
simulation. All of these observables require finding intersections
of the walk with a given curve (a line, parabola or circle). 
Searching through the walk one step at a time for these intersections
would be disastrous, since it would require a time $O(N)$. Such a
search is easily avoided. At a given site in the walk we do not simply
check if the next step intersects the curve. Instead we compute the 
distance from the site to the curve. The walk must take at least this
many steps before it can intersect the curve, so we can jump ahead this
many steps in the walk before we check again for an intersection.

\begin{appendix}
\section{Proof of irreducibility}

In this appendix we prove that the pivot algorithm is irreducible in the
half-plane and cut-plane that we have been considering. The proof is 
very similar to the proof for the full plane \cite{ms}.
We show that for any self-avoiding walk in the half-plane 
(cut-plane, respectively), there is a sequence of pivot operations  
which ``unfold'' the walk into a straight line and such that each 
walk produced in this unfolding process is self-avoiding and remains 
in the half-plane (cut-plane, respectively). 

We first consider the half-plane. The restriction is that except for the 
starting point of the walk at the origin, the walk must remain strictly 
above the horizontal axis. We will show that the number of 
turns in the walk can be decreased by one. 
We denote the sites in the walk by $\omega(i)$ where $i=0,1,\cdots,N$. 
We will say there is a turn at $\omega(i)$ if $\omega(i-1), \omega(i)$ 
and $\omega(i+1)$ are not co-linear.  

We will consider cases based on the direction of the last step of the walk.
If it is to the right, i.e., $\omega(N)=\omega(N-1)+(1,0)$ we proceeds
as follows. Let $l$ be the largest integer such that the line 
$y-x=l$ contains a site on the walk. So the walk is entirely below
or on this line. Let $i$ be the largest integer such that $\omega(i)$ 
is on this line. Since the last step of the walk is to the right, 
$\omega(N)$ is not on this line. So $i<N$. 
Since the first step of the walk in the half-plane must be up, $i$ 
cannot be $0$. 
We take $\omega(i)$ as the pivot point and reflect the portion of the walk 
from $\omega(i)$ to $\omega(N)$ in the line $y-x=l$. The reflected
portion of the walk lies entirely above the line, so the reflection
does not produce self-intersections. Furthermore, since the walk 
was on or below the line, the reflection can only increase the $y$ coordinates
of points on the walk. So the new walk is still in the upper half plane. 
The walk before this reflection
has a turn at $\omega(i)$ and the reflected walk does not. The reflection
does not add any turns to the walk, so the total number of turns 
decreases by one. If the final step of the walk is to the left, 
i.e., $\omega(N)=\omega(N-1)-(1,0)$, we use an analogous procedure 
with lines $y+x=l$ to reduce the number of turns in the walk. 

Now suppose that the final step of the walk is either up or down, i.e., 
$\omega(N)=\omega(N-1) \pm (0,1)$. Consider the vertical line which 
contains this last step. First suppose that the walk lies entirely 
to the right of or on this vertical line. Let $i<N$ be the largest integer
such that there is a turn at $\omega(i)$. (Of course, if there are 
no turns the walk is a straight line and we are done.)
The walk is a straight segment 
from $\omega(i)$ to $\omega(N)$ which lies on the vertical line. 
We take $\omega(i)$ as the pivot point and perform a rotation of 
$90$ degrees ($-90$, respectively) if the last step of the walk is 
up (down, respectively). This rotates the segment from $\omega(i)$ to 
$\omega(N)$ to the left of the vertical line and removes the turn at 
$\omega(i)$. No new turns are added to the walk, so the total number of 
turns decreases by one. If the walk likes entirely to the right or on 
the vertical line containing the last step, an analogous argument 
shows the number of turns can be reduced by one. 

Now suppose that the walk contains sites on both sides of the vertical 
line which contains the last step of the walk. Let $d$ be the horizontal 
width of the walk:
\be 
d = \max \{ x : (x,y) = \omega(i), {\rm \, for \, \, some \, \, } \, i, y \}
- \min \{ x : (x,y) = \omega(i), {\rm \, for \, \, some \, \, } \, i, y \}
\ee 
We will show that $d$ can be increased. 
Let $l$ be the smallest integer 
such that the vertical line $x=l$ contains sites in the walk. So the walk 
lies on or to the right of this line. Note that
$\omega(N)$ is not on this line. Let $i<N$ be the largest integer such that 
$\omega(i)$  is on this line. We take $\omega(i)$ as the pivot point and 
reflect the walk from $\omega(i)$ to $\omega(N)$ in the line $x=l$. 
This increases the width of the walk. (The argument is the same as that 
given in \cite{ms}.) The reflection does not change the y-coordinate 
of points on the walk, so the new walk is still in the upper half-plane.
Note that in the new walk the last step is in the 
same direction as before, i.e., either up or down. 
So we can repeat this procedure to increase $d$ until we obtain 
a walk which lies entirely on or
to one side of the vertical line containing the final step. When we reach
such a walk we apply the procedure of the proceeding paragraph to 
reduce the number of turns by one. This completes the proof for the 
case of the half-plane. 

Now consider the cut-plane. The restriction now is that the walk cannot 
contain sites of the form $(x,0)$ with $x \ge 0$, except for the starting 
point at the origin. 
We again consider cases based on the direction
of the last step of the walk. If it is to the right, we proceed as in the 
half-plane algorithm. Note that the line involved, $y-x=l$, must have 
$l \ge 0$ since the walk starts at the origin. (In fact $l$ must be at least
$1$, but we do not need this.) The reflected portion of the walk will lie
above this line while the cut, the non-negative real axis, lies below it.
So the reflection produces a walk that lies in the cut-plane. 

If the last step of the walk is to the left, a different procedure is 
needed to avoid producing a walk that intersects the cut. 
Consider the lines $x-y=l$ and $x+y=l$. They intersect at $(l,0)$ and 
divide the plane into four quadrants which we will describe as being
left, right, above and below the point $(l,0)$. We take $l$ to be the 
smallest integer such that the sites on the walk lie 
in the quadrant to the left of $(l,0)$ or on the lines. 
($l$ is necessarily positive.) 
We then let $i$ be the largest integer such that $\omega(i)$ is on one of the 
two lines. (It is not $N$ since the last step of the walk is to the left.)
Note that $\omega(i)$ cannot be $(l,0)$. 
We take $\omega(i)$ as the pivot point and reflect the walk from 
$\omega(i)$ to $\omega(N)$ in the line containing $\omega(i)$. 
The reflected portion of the walk will lie either in the quadrant above
or below $(l,0)$ and so cannot intersect the cut. 

If the last step of the walk goes up or down we use the algorithm for 
the half-plane. There is a subtle point here. Recall than when   
the walk has points on both sides of the vertical line containing the last
step, we chose $l$ so that the walk is to the right of or on the vertical 
line $x=l$. For the half plane we could have chosen it so that ``right'' is 
replaced by ``left.''  
For the cut-plane this choice could result in a reflected
walk that intersects the cut. To see that our choice does not produce 
a walk that intersects the cut, we observe that  
since the walk starts at the origin, $l$ must be negative. 
The reflected portion of the walk will lie on or to the left of the line 
$x=l$, and so will not intersect the cut. 
This completes the proof for the cut-plane. 

\end{appendix}

\bigskip
\bigskip

\noindent {\bf Acknowledgments:}
Oded Schramm pointed out to me that one can find the conformal map 
from the half-plane with a circular arc removed onto the half-plane
and so compute the distribution of $\Theta_e$. 
The author thanks the Mathematics Department of 
Universit\'e Paris-Sud (Orsay) and the Institute des Hautes 
\'Etudes Scientifiques (Bures-sur-Yvette), 
where some of this work was performed, for their hospitality.
This work was supported by the National Science Foundation (DMS-9970608 and 
DMS-0201566).

\end{document}